\documentclass[11pt]{article}
\usepackage{amssymb,amsmath}
\usepackage[latin1]{inputenc}
\DeclareMathOperator{\op}{op}

\begin{document}
\newcommand{\hhom}{{\mathbb H}}
\newcommand{\R}{{\mathbb R}}
\newcommand{\N}{{\mathbb N}}
\newcommand{\Z}{{\mathbb Z}}
\newcommand{\C}{{\mathbb C}}
\newcommand{\T}{{\mathbb T}}
\newcommand{\rn}{{\mathbb R}^n}

\newcommand{\cA}{{\mathcal A}}
\newcommand{\cB}{{\mathcal B}}
\newcommand{\cC}{{\mathcal C}}
\newcommand{\cD}{{\mathcal D}}
\newcommand{\cE}{{\mathcal E}}
\newcommand{\cF}{{\mathcal F}}
\newcommand{\cG}{{\mathcal G}}
\newcommand{\cH}{{\mathcal H}}
\newcommand{\cI}{{\mathcal I}}
\newcommand{\cJ}{{\mathcal J}}
\newcommand{\cK}{{\mathcal K}}
\newcommand{\cN}{{\mathcal N}}
\newcommand{\cL}{{\mathcal L}}
\newcommand{\cP}{{\mathcal P}}
\newcommand{\cQ}{{\mathcal Q}}
\newcommand{\cS}{{\mathcal S}}
\newcommand{\cT}{{\mathcal T}}

\newcommand{\ga}{\alpha}
\newcommand{\gb}{\beta}
\renewcommand{\gg}{\gamma}
\newcommand{\gG}{\Gamma}
\newcommand{\gd}{\delta}
\newcommand{\eps}{\varepsilon}
\newcommand{\gve}{\varepsilon}
\newcommand{\gk}{\kappa}
\newcommand{\gl}{\lambda}
\newcommand{\gL}{\Lambda}
\newcommand{\go}{\omega}
\newcommand{\gO}{\Omega}
\newcommand{\gvp}{\varphi}
\newcommand{\gt}{\theta}
\newcommand{\gT}{\Theta}
\renewcommand{\th}{\vartheta}
\newcommand{\gs}{\sigma}

\newcommand{\fT}{{\mathfrak T}}

\newcommand{\ol}{\overline}
\newcommand{\ul}{\underline}

\newcommand{\Dbar}{D\hspace{-1.5ex}/\hspace{.4ex}}

\newcommand{\Pf}{{\em Proof.}~}
\newcommand{\Proof}{\noindent{\em Proof}}
\newcommand{\pitensor}{\hat\otimes_\pi}
\newcommand{\eproof}{{~\hfill$ \triangleleft$}}
\renewcommand{\i}{\infty}
\newcommand{\rand}[1]{\marginpar{\small  #1}}
\newcommand{\forget}[1]{}

\newcommand{\cut}{C^\infty_{tc}(T^-X)}
\newcommand{\Ctc}{C^\infty_{c}}
\newtheorem{leer}{\hspace*{-.3em}}[section]
\newenvironment{rem}[2]%
{\begin{leer} \label{#1} {\bf Remark. }  
 {\rm #2 } \end{leer}}{}
\newenvironment{lemma}[2]%
{\begin{leer} \label{#1} {\bf Lemma. }  
{\sl #2} \end{leer}}{}
\newenvironment{thm}[2]%
{\begin{leer}\label{#1} {\bf Theorem. } 
{\sl #2} \end{leer}}{}
\newenvironment{dfn}[2]%
{\begin{leer}  \label{#1} {\bf Definition. }
{\rm #2 } \end{leer}}{}
\newenvironment{cor}[2]%
{\begin{leer} \label{#1} {\bf Corollary. }
{\sl #2 } \end{leer}}{}
\newenvironment{prop}[2]%
{\begin{leer} \label{#1} {\bf Proposition. }
{\sl #2} \end{leer}}{}
\newenvironment{extra}[3]%
{\begin{leer} \label{#1} {\bf #2. } 
{\rm #3 } \end{leer}}{}
\renewcommand{\theequation}{\thesection.\arabic{equation}}
\renewcommand{\labelenumi}{{\rm (\roman{enumi})}}

\newcounter{num}
\newcommand{\bli}[1]{\begin{list}{{\rm(#1{num})}\hfill}{\usecounter{num}\labelwidth1cm
\leftmargin1cm\labelsep0cm\rightmargin1pt\parsep0.5ex plus0.2ex minus0.1ex
\itemsep0ex plus0.2ex\itemindent0cm}}
\newcommand{\eli}{\end{list}}

\def\Im{{\rm Im}\,}
\def\lra{\longrightarrow}
\def\Re{{\rm Re}\,}
\def\rpb{\overline\R_+}
\def\sumj#1{\sum_{j=0}^{#1}}
\def\sumk#1{\sum_{k=0}^{#1}}
\def\vect#1#2#3{\begin{array}{c}#1\\#2\\#3\end{array}}
\def\vec#1#2{\begin{array}{c}#1\\#2\end{array}}
\def\skp#1{\langle#1\rangle}

\title{A Continuous Field of $C^*$-algebras and the Tangent Groupoid for  Manifolds with Boundary}
\author{\sc Johannes Aastrup, Ryszard Nest and Elmar Schrohe}
\date{}
\maketitle 

{\small
{\bf Abstract.} For a smooth manifold $X$ with boundary we construct a semigroupoid $\cT^-X$ 
and a continuous field $C^*_r(\cT^-X)$ of $C^*$-algebras which extend Connes' construction
of the tangent groupoid.
 
We show the asymptotic multiplicativity of $\hbar$-scaled truncated pseudodifferential operators 
with smoothing symbols and compute the $K$-theory of the associated symbol algebra. 

{\bf Math.\ Subject Classification} 58J32, 58H05, 35S15, 46L80.

}
\tableofcontents

\section*{Introduction}
It is a central idea of semi-classical analysis to consider Planck's constant $\hbar$ as a small real variable 
and to study the relation between systems in mechanics and systems in quantum mechanics by associating to 
a function $f=f(x,\xi)$ on the cotangent bundle of a manifold 
the $\hbar$-scaled pseudodifferential operator $\op_\hbar(f)$ with symbol $f(x,\hbar\xi)$ and 
analyzing their relation as $\hbar\to0$.   

For  $f \in \cS (T^*\R^n)$, for example,  a basic  estimate states that
\begin{equation}\label{0.1}
\lim_{\hbar \to 0}\| \op_\hbar (f )\|=\| f \|_{\rm sup} .
\end{equation} 
Moreover, given a second symbol $g\in \cS (T^*\R^n)$ we have 
\begin{equation}\label{0.2}
\lim_{\hbar \to 0}\| \op_\hbar (f)\op_\hbar (g)-\op_\hbar (fg)\|=0; 
\end{equation} 
in other words, the map $\op_\hbar $ is asymptotically multiplicative.

As both statements concern the asymptotic behavior of pseudodifferential operators, 
it is somewhat surprising that they can be proven within the framework 
of continuous fields of $C^*$-algebras associated to amenable Lie groupoids, more precisely, the $C^*$-algebra of 
the so-called tangent groupoid $\cT M$,
cf.\ Connes \cite[Section II.5]{Connes94}.

For a boundaryless manifold $M$,  $\cT M $ is constructed
by gluing the tangent space $TM$ to the cartesian product 
$M\times M\times ]0,1]$ via the map $TM\times [0,1]\ni (m,v,\hbar)\mapsto (m,\exp_m(-\hbar v),\hbar)$. 
It has the natural cross-sections $\cT M(\hbar)$, $0\le \hbar\le 1$, 
given by $TM$ for $\hbar =0$ and by $M\times M\times \{\hbar\}$ for $\hbar\not=0$.

The basic observation, establishing the link between 
$\hbar$-scaled pseudodifferential operators and
the tangent groupoid, is the following:
In the Fourier transformed picture, 
the $\hbar$-scaled pseudodifferential operator $\op_\hbar(f)$ 
becomes the convolution operator $\rho_\hbar(\hat f)$ acting by
$$\rho_\hbar(\hat f)\xi(x)
=\frac1{\hbar^n}\int \hat f\big(x,\frac{x-y}\hbar\big)\xi(y)dy,\quad\xi\in L^2(\R^n), $$
and for $\hbar \not=0$, the mappings $\rho_\hbar$ 
(or better their generalization to the manifold case) 
coincide with the natural representations of 
$C^\infty_c(\cT M(\hbar))$ by convolution operators. 

The $\rho_\hbar$, $\hbar\not=0$, are complemented by the representation $\pi_0$ of $C^\infty_c(TM)$ on
$L^2(TM)$ via convolution in the fiber 
which in turn coincides with the natural representation of $C^\infty_c( \cT M(0)).$ 

Now the tangent groupoid is additionally  amenable, so that,
according to a theorem by Anantharaman-Delaroche and Renault \cite{AnDeRe},
the reduced  $C^*$-algebra $C^*_r(\cT M)$, defined as the closure of $C^\infty_c(\cT M)$ with respect to the 
natural representations,    
and the full $C^*$-algebra $C^*(\cT M)$, i.e., the 
closure with respect to all involutive Hilbert space representations, are isomorphic.

The crucial fact then is that $C_r^*(\cT M)$ is a continuous field of $C^*$-algebras over $[0,1]$; the fiber over $\hbar$ is
$C^*_r(\cT M(\hbar))$. 
An elegant way to establish the continuity is to 
show upper semi-continuity and lower semi-continuity separately, 
noticing that upper semi-continuity is easily proven in $C^*(\cT M)$ while lower semi-continuity 
is not difficult to show in $C_r^*(\cT M)$. As both $C^*$-algebras are isomorphic, continuity follows. 
For a good account of these facts see  
\cite{LandsmanRamazan} by Landsman and Ramazan.
The identities \eqref{0.1} and \eqref{0.2} are then an immediate consequence
of the continuity of the field.

In the present paper we consider manifolds with boundary. 
The analog of the usual pseudodifferential calculus  here
is Boutet de Monvel's calculus for boundary value problems \cite{MR53:11674}. 
In order to obtain an operator algebra, one cannot work with pseudodifferential 
operators alone, but has to introduce an additional class of operators, the so-called singular Green
operators. The reason is the way pseudodifferential operators act on functions defined on a half space:
One first extends the function (by zero) to the full space, 
then applies the pseudodifferential operator and finally restricts
the result to the half space again -- one often speaks of trucated pseudodifferential operators. 
Given two pseudodifferential operators $P$ and $Q$, the `leftover operator' $L(P,Q)=(PQ)_+-P_+Q_+$, i.e.\ the 
difference between
the trucated pseudodifferential operator $(PQ)_+$ associated to the composition $PQ$ and the composition of the   
truncated operators $P_+$ and $Q_+$ associated with $P$ and $Q$ is a typical example of such a singular Green operator. 
The singular Green operators `live' at the boundary. They are smoothing operators in the interior, while, close to 
the boundary, they can be viewed as operator-valued pseudodifferential operators along the boundary, acting like smoothing
operators in the normal direction.

In the full algebra which consists -- at least in the slightly simplified picture we have here -- of sums of (truncated)
pseudodifferential operators and singular Green operators, the singular Green operators form an ideal. 

With this picture in our mind, we construct an analog of Connes' tangent groupoid for a manifold $X$ with boundary.
Our semi-groupoid $\cT^- X$ consists of the groupoid $X\times X\times ]0,1]$ to which we glue, with the same map as above, the 
half-tangent space $T^-X$, which comprises all those tangent vectors to $X$ for which $\exp_m(-\hbar v)$ lies in $X$ for
small $\hbar$ (note that this condition is only effective at the boundary of $X$). 
As before, we have natural cross-sections $\cT^- X(\hbar)$, coinciding with $X\times X\times\{\hbar\}$ for $\hbar\not=0$ and 
with $T^-X$ for $\hbar=0$. 

For $\hbar\not=0$, the operators $\rho_\hbar$ (with integration now restricted to $X$), are the 
natural representations of the groupoid $\cT^-X(\hbar)$. 
At $\hbar=0$ we use two mappings. 
The first, $\pi_0$ is the analog of the above map $\pi_0$. It acts on the 
tangent space of $X$ by convolution. 
The second, $\pi_0^\partial$, acts on the half tangent space over the boundary by half-convolution: 
$\pi_0^\partial:C^\infty_c(T^-X)\to\cL(L^2(T^-X|_{\partial X})$ is given by 
$$\pi_0^\partial(f)\xi(m,v)=\int _{T^-_mX}f(m,v-w)\xi(m,w)dw.$$

In order to avoid problems concerning the topology of $\cT^-X$, 
we denote by $C^\infty_c(\cT^-X)$ the space
of all restrictions of functions in $C_c^\infty(\cT\widetilde X)$ to $\cT^-X$; here $\widetilde X$ is 
a boundaryless manifold containing $X$.

The reduced $C^*$-algebra $C^*_r(\cT^-X)$ is then defined as the $C^*$-closure of $C^\infty_c(\cT^-X)$ with 
respect to the $\rho_\hbar$, $\hbar\not=0$, and $\pi_0,\pi_0^\partial$ for $\hbar=0$. For the full $C^*$-algebra 
we use all involutive representations.

We show that $C^*_r(\cT^-X)$ is a continuous field of $C^*$-algebras over $[0,1]$, where the fiber over
$\hbar\not=0$ is $C^*_r(\cT^-X(\hbar))$, and the fiber over $\hbar=0$ is the $C^*$-closure of 
$C^\infty_c(T^-X)$ with respect to  $\pi_0$ and $\pi_0^\partial$.

The proof of continuity is again split up into showing upper semi-continuity and lower semi-continuity.
According to an idea by Rieffel \cite{Rieffel}, lower  semi-continuity is established using 
strongly continuous representations. 
The basic idea for the proof of upper semi-continuity would be to infer an isomorphism between $C^*_r(\cT^-X)$ and
$C^*(\cT^-X)$ from the amenability of $\cT^-X$. 
However, as $T^-X$ is only a semi-groupoid, we make a little detour: 
Using short exact sequences and the amenability of the tangent groupoids for boundaryless manifolds 
we prove that $C^*_r(T^-X)$ is isomorphic to the closure of $C^\infty_c(T^-X)$ with respect to the 
upper semi-continuous norm.

The present study should be seen as a step towards fitting Boutet de Monvel's calculus for boundary value problems 
into the framework of deformation quantization and groupoids, 
in the spirit of Connes \cite{Connes94}, Monthubert and Pierrot \cite{MonthuberPierrot},
Nest and Tsygan \cite{MR1350407}, 
\cite{MR1337107}, Nistor and Weinstein and Xu \cite{MR1687747}, 
Eventually one could hope to develop an algebraic index theory for these deformations in the spirit of Nest and Tsygan. 

The structure of the paper is as follows:
In the first section we review the case of boundaryless manifolds. 
We introduce the basic notions and show how \eqref{0.1} and \eqref{0.2} are derived 
with the help of  the continuous field of $C^*$-algebras associated to the tangent groupoid.

We then consider a manifold $X$ with boundary. In order to make the presentation more transparent, we
first study the case where $X=\R^n_+=\{(x_1,\ldots,x_n)\ |\ x_n\ge0\}$. 
Here all relevant features show up, but computations are easier to perform. 
We then go over to the general case.

In Section 3 we determine the $K$-theory of the symbol algebra $C^*_r(T^-X)$.
Starting from the short exact sequence 
$$0\longrightarrow C^*_r(TX^\circ)\longrightarrow C^*_r(T^-X)\longrightarrow Q\longrightarrow 0$$
we show that the quotient $Q$ can be identified with $C_0(T^*\partial X)\otimes \cT_0$, where 
$\cT_0$ is an ideal in the Toeplitz algebra with vanishing $K$-theory. In particular, we obtain the isomorphism 
$$K_i(C^*_r(T^-X))\cong K_i(C_0(T^*X)),\quad i=0,1.$$

The appearance of the Toeplitz operators can be seen as a feature inherent in the geometry of the problem. 
In fact, the construction of an algebra of pseudodifferential operators on a closed (Riemannian) manifold 
amounts to the construction of a suitably completed operator algebra, generated by 
multivariable functions of vector fields and the operators of multiplication by smooth functions.

In the boundaryless case, one can localize to $\R^n$ and reduce the task essentially to defining 
$f(D)$ for a classical symbol $f$ and $D=(D_1,\ldots, D_n)$ 
with the  vector fields $D_j=i\partial_{x_j}$. 
One convenient way of achieving this is to use the operator families $e^{itD_j}$ and to  let 
$$
f(D)=(2\pi)^{-n}\int \widehat f(\xi)e^{i\xi D}\,d\xi
$$ 
with the Fourier transform $\widehat f$  of $f$ and $\xi D=\xi_1D_1+\ldots +\xi_n D_n$. 
Note that the use of the $e^{i\xi D}$ is purely geometric and only relies on the fact that 
vector fields integrate to flows. 

On a manifold with boundary, one will have vector fields transversal 
to the boundary which do not integrate to flows. In this case, one has two possibilities: 
The first is to restrict the class of admissible vector fields to those which {\em} do 
integrate. 
This is a basic idea in the pseudodifferential calculi introduced by Melrose \cite{MelroseKyoto}, 
see also Ammann, Lauter, Nistor \cite{AmmannLauterNistor}.    

In Boutet de Monvel's calculus, on the other hand, transversal vector fields are admitted. 
After localization to $\ol \R^n_+$, we may focus on $D_n$.
One of the functions one would certainly like to define is the Cayley transform 
(recall that the Cayley transform $C(A)$ of an operator $A$ is given by 
$C(A)= (A-i)(A+i)^{-1}=1-2i(A+i)^{-1}$). 

Now it is well-known that the Cayley transform $C(A)$ is an isometry, 
and that it is a unitary if and only if  $A$ is selfadjoint. 
As there is no selfadjoint extension of $D_n$,
its Cayley transform will be a proper isometry.
Hence by a theorem of Coburn \cite{Coburn,Coburn2}, 
the algebra generated by it (which becomes part of the calculus),
is the Toeplitz algebra. 

While the pseudodifferential calculus for closed manifolds is commutative modulo lower order
terms, this calculus is not. 
From a geometric point of view, the resulting algebra can thus be seen as a noncommutative   
completion of the manifold with boundary.

{\em Remark on the notation.} A variety a representations naturally comes up in this 
context. In order to distinguish their origin, we will apply the following rule.  
Representations related to the groupoid structure are denoted by $\pi$ (possibly indexed), 
asymptotic pseudodifferential operators by $\rho_\hbar$ and the asymptotic 
Green operators (introduced in Section 2) by $\kappa_\hbar$. 

\section{The Classical Case}
\subsection*{Groupoids}
A groupoid $G$ is a small category where all the morphisms are invertible. 
We will denote by $G^{(0)}$ the set of objects in $G$ and by $G^{(1)}$ the set of morphisms. 
We will also call $G^{(0)}$ the base and the elements in $G^{(1)}$ the arrows. 
On $G^{(1)}$ there are two maps $r,s$ into $G^{(0)}$. The first map, $r $,  is the range object of a morphism 
and the second, $s$, the source. 
For $x\in G^{(0)}$ we define $G^x=r^{-1}(x)$ and $G_x=s^{-1}(x)$. 
There is an embedding $\iota$ of $G^{(0)}$ into $G^{(1)}$ given by mapping an object 
to the identity morphism on this object. 
Furthermore we define $G^{(2)}$ to be the subset of composable morphisms of $G^{(1)}\times G^{(1)}$.

\begin{dfn}{liegr}{A Lie groupoid $G$ is a groupoid together with a manifold structure on $G^{(0)}$ and $G^{(1)}$ such that the maps $r,s$ are submersions,  the map $\iota$ and the the composition map $G^{(2)} \rightarrow G^{(1)}$ are smooth. }
\end{dfn}  
To a given a smooth  manifold $M$ without boundary there are associated two canonical Lie groupoids. 
The first is the tangent bundle $TM$ of $M$. The groupoid structure is given by 
\begin{eqnarray*}
G^{(0)} =M, && G^{(1)}=TM  \\
r(m,X)=m, && s(m,X)=m\\
(m,X)\circ (m,Y)&=&(m,X+Y)
\end{eqnarray*}
The second one is the pair groupoid $M\times M$ with
\begin{eqnarray*}
G^{(0)}=M, && G^{(1)}=M\times M \\
r(m_1,m_2)=m_1, &&s(m_1,m_2)=m_2 \\
(m_1,m_2)\circ (m_2,m_3)&=&(m_1 , m_3)
\end{eqnarray*}
Both  are clearly Lie groupoids. 

\extra{Haar}{Haar systems}{A smooth left Haar system on a Lie groupoid is a family of measures 
$\{ \lambda^x \}_{ x\in G^{(0)}}$ on $G$ with $\hbox{supp}\lambda^x=G^x$ which is left invariant, i.e. 
$\gamma(\lambda^{s(\gamma )})=\lambda^{r(\gamma )}$, and for each $  f\in C^\infty_c (G^{(1)})$, 
the function on $G^{(0)}$ defined by
$$x \mapsto \int fd\lambda^x, \quad  f\in C^\infty_c (G^{(1)})$$ 
is smooth.
In \cite[Proposition 3.4]{LandsmanRamazan}, it is proven that all Lie groupoids possess a smooth left Haar system.
Similarly, a right Haar system $\{\gl_x\}$ is given by $\gl_x=(\gl^x)^{-1}$.
} 

\begin{dfn}{ame}{A   Lie groupoid $G$ with a smooth left Haar system $\lambda^x$ 
is called topologically amenable if there exists a net of nonnegative continuous functions 
$\{f_i\}$ on $G^{(1)}$ such that
\begin{enumerate}
\item For all $i$ and for $x\in G^{(0)}$, $\int f_i d\lambda^x=1$.
\item The functions $\gamma \mapsto \int | f_i(\gamma^{-1}\gamma')-f_i(\gamma')|d\lambda^{r(\gamma)}(\gamma')$ converge 
uniformly to zero on compact subsets of $G^{(1)}$.
\end{enumerate}}
\end{dfn}
It is easy to verify that the two groupoids $TM$ and $M\times M$ are topologically amenable.

\extra{tg}{Connes' tangent groupoid}{Let $M$ be a smooth manifold. 
Connes tangent groupoid $\cT M$ is a blow up of the diagonal in $M\times M$. More specifically:

Let $\cT M=TM \cup (M\times M\times ]0,1])$ as a set. 
The groupoid structure is just the fiberwise groupoid structure 
coming from the groupoid structure on $TM$ and $M\times M$. 
The manifold structure on $M\times M \times ]0,1]$ is obvious. 
We next glue $TM$  to $M\times M \times ]0,1]$ to get a manifold structure on $\cT M$. 
To this end we choose a Riemannian metric  on $M$ and glue with the charts
$$TM\times [0,1]\supseteq U\ni (m,v,\hbar)\mapsto \left\{ 
\begin{array}{cc}
(m,v) & \hbox{for }\hbar =0 \\
(m, \exp_m(-\hbar v),\hbar)& \hbox{for }\hbar \not= 0,
\end{array}
\right.$$
where $\exp_m$ denotes the exponential map and $U$ is an open neighborhood of $M\times \{0\}\subset TM\times \{0\}$; 
here, $M$ is embedded as the zero section.

Here, $G^{(0)}=M\times [0,1]$. For $\hbar\not=0$ and  $x=(\tilde m,\tilde \hbar)\in G^{(0)}$,   
we have $G^x=\{(\tilde m,m,\tilde \hbar): m\in M\}$; for $x=(\tilde m,0)$, $G^x=T_{\tilde m}M$. 
Fixing  the measure $\mu$ on $M$ induced by the metric, 
we obtain a Haar system  $\{\lambda^x\}_{x\in G^{(0)}}$ by $\lambda^{(\tilde m,\tilde \hbar)}
=\hbar^{-n}\mu$, $\hbar \not=0$; for $\hbar=0$, 
we let $\lambda^{(\tilde m,0)}$  be the measure on $T_mM$ given by the metric.

This makes $\cT M$ a Lie groupoid, see  \cite{LandsmanRamazan}. 
}

\extra{c*}{C*-algebras associated to groupoids}{Let $G$ 
be a Lie groupoid with a smooth left Haar system $\lambda$. 
On $C^\infty_c(G^{(1)})$ we define a $*$-algebra structure by
\begin{eqnarray}
(f*g)(\gamma)&=&\int_{G^{s(\gamma)}}f(\gamma \gamma_1)g(\gamma_1^{-1})\ d\lambda^{s(\gamma )}(\gamma_1) \\
f^*(\gamma )&=&\overline{f(\gamma^{-1})}
\end{eqnarray} 
There are involutive representations $\pi_x$, $x\in G^{(0)}$,  
of this $*$-algebra on the Hilbert spaces $L^2(G_x,\lambda_x)$ given by
\begin{equation}\label{repsl}
\pi_x(f)\xi(\gamma )=\int_{G^{x}}f(\gg\gamma_1)\xi (\gamma_1^{-1})~d\lambda^{x}(\gamma_1),
\qquad \xi \in L^2(G_x,\lambda_x).
\end{equation}
}

\begin{dfn}{deffcg}{The full $C^*$-algebra $C^*(G)$ of a groupoid is the $C^*$-completion 
of the $*$-algebra $C^\infty_c(G^{(1)})$ with respect to all involutive Hilbert space representations.

The reduced $C^*$-algebra $C^*_r(G)$ of $G$ is the $C^*$-completion of $C_c^\infty(G)$ 
with respect to the representations \eqref{repsl}.}  
\end{dfn}

Note that, by universality, we have a quotient map from $C^*(G)$ to $C^*_r(G)$. 

\begin{rem}{1.7}{Although the construction of the $*$-algebra 
structure on $C_c^\infty (G^{(1)})$ and the representations \eqref{repsl}  
use a smooth Haar system, the algebra is independent of the choice. 
See \cite{LandsmanQTG} for a detailed exposition.}
\end{rem}  

\extra{tanfour}{Example}{For the tangent bundle $TM$ of a manifold, 
the space $G_m$ is just $T_mM$ and the representation is
$$\pi_m (f)\xi(v)=\int_{T_mM} f(m,v-w)\xi (w)\, dw\quad \xi \in L^2(T_mM).$$
By Fourier transform in each fiber $T_mM$, the $C^*$-algebra $C^*_r(TM)$ 
becomes isomorphic to $C_0(T^*M)$, the continuous functions on $T^*M$ 
vanishing at infinity. }

The importance of topological amenability lies in the following result from \cite{AnDeRe}:
\begin{prop}{amec}{When $G$ is topologically amenable
 the quotient map from $C^*(G)$ to $C^*_r(G)$ is an isomorphism.}
\end{prop}

\subsection*{Continuous Fields and $\hbar$-Scaled Pseudodifferential Operators} 

\dfn{cf}{A continuous field of $C^*$-algebras 
$(A,\{A(\hbar),\gvp_\hbar\}_{\hbar\in[0,1]})$ over $[0,1]$ consists of a $C^*$-algebra $A$, $C^*$-algebras
$A(\hbar)$, $\hbar\in[0,1]$, with surjective homomorphisms $\gvp_\hbar:A\to A(\hbar)$ 
and an action of $C([0,1])$ on $A$ such that for all $a\in A$
\begin{enumerate}
\item 
The function $\hbar\mapsto \|\gvp_\hbar(a)\|$ is continuous;
\item 
$\|a\|=\sup_{\hbar\in[0,1]}\|\gvp_\hbar(a)\|$;
\item 
For $f\in C([0,1])$, $\gvp_\hbar(fa)=f(\hbar)\gvp_\hbar(a)$.
\end{enumerate}}

\thm{cfprop}{For the tangent groupoid $\cT M$ we define 
$\cT M(0)=TM$ and $\cT M(\hbar)=M\times M\times \{\hbar\}$ for $\hbar\not=0$.
 The pullback
under the inclusion $\cT M(\hbar)\hookrightarrow \cT M$ induces a map  
$\gvp_{\hbar}:C^\infty_c(\cT M)\to C^\infty_c(\cT M(\hbar))$ which
 extends by continuity to a surjective $*$-homomorphism $\gvp_{\hbar}:
C^*_r(\cT M)\to C^*_r(\cT M(\hbar))$. 
The $C^*$-algebras $A= C^*_r(\cT M)$ and $A(\hbar)=C^*_r(\cT M(\hbar))$ with the maps $\gvp_{\hbar}$
form a continuous field over $\R$.}

\Proof. Together with the amenability of $\cT M$ and Proposition \ref{amec} 
this is immediate from Theorem 6.4 in \cite{LandsmanRamazan}.\eproof

\extra{1.12}{$\hbar$-scaled pseudodifferential operators}{For $0<\hbar\le 1$ define 
 $\rho_\hbar: C^\infty_c(T\R^n)\to \cL(L^2(\R^n))$  by
\begin{eqnarray}
\rho_\hbar (f)\xi(x) &=& \int f(x,w)\xi(x-\hbar w)\,dw
=\hbar^{-n}\int f\big(x,\frac{x-w}\hbar\big)\xi(w)\,dw, \quad \xi\in L^2(\R^n)
\end{eqnarray}
We complement this by the map $\pi_0: C^\infty_c(T\R^n)\to \cL(L^2(T\R^n))$
\begin{eqnarray}\label{tilderho}
 \pi_0(f)\xi(x,v)&=& \int f(x,w)\xi(x,v- w) dw.
\end{eqnarray}
}

\rem{1.12a}{(a) We can define $\tilde\rho_\hbar: C^\infty_c(T\R^n)\to \cL(L^2(T\R^n))$, $\hbar\ge 0$, by 
$$\tilde \rho_\hbar(f)\xi(x,v)= \int f(x,w)\xi(x-\hbar w,v- w) dw$$
and then obtain a more consistent representation. Note that
for $h>0$ the representations $\rho_\hbar$ and $\tilde \rho_\hbar$ are unitarily equivalent. 

(b)
On a smooth Riemannian manifold $M$ we define $\rho_\hbar$ by  
\begin{eqnarray}
(\rho_\hbar f)\xi(x) &=& \int \psi (x,\exp_x(-\hbar w)) f(x,w)\xi(\exp_x(-\hbar w))dw\nonumber\\
&=&\hbar^{-n}\int \psi(x,y) f(x,-\exp^{-1}(x,y)/\hbar)\xi(y)\, dy.\label{eq1.8}
\end{eqnarray}
Here  $\psi \in C^\infty (M\times M)$ is a function,
which is one on a neighborhood of the diagonal, $0\leq \psi \leq 1$ and such that 
$$\exp :TM \rightarrow M\times M,\quad (m,v)\mapsto (m,\exp_mv),$$ 
maps a neighborhood of the zero section diffeomorphically to the support of 
$\psi$; a similar construction applies to $\tilde \rho$.

Note that for two representations $\rho^1_\hbar,\rho^2_\hbar$, defined with cut-off functions $\psi_1$ and $\psi_2$, the norm $\| \rho_\hbar^1 (f)-\rho_\hbar^2(f)\|$ tends to zero as $\hbar \to 0$.}

\lemma{1.13}{To each $f\in C^\infty_c(TM)$ we associate a function $\tilde f\in C^\infty(\cT M)$ by
$$
\begin{array}{ll}\tilde f(x,v,0)=f(x,v)&\text{ for  } \hbar=0,x\in M,v\in T_xM;\\
\tilde f(x,y,\hbar)=\psi (x,y) f(x,-\exp^{-1} (x,y)/ \hbar)& \text{ for } \hbar\not=0,x,y\in M.
\end{array}$$ 
By \eqref{eq1.8}  
$$
\pi_{(x,\hbar)}(\tilde f)=\rho_{\hbar}(f) \text{ and } 
\|\gvp_{\hbar}(\tilde f)\|_{\cT M(\hbar)}=\sup_x\|\pi_{(x,\hbar)}(\tilde f)\|_{\cL (L^2(G_{(x,\hbar)},\lambda_{(x,\hbar)}))}
=\|\rho_\hbar(f)\|.
$$
}

\thm{1.14}{We denote by $\widehat f$ the Fourier transform of $f$ with respect to the covariable.
Then \\
{\rm (a)} $\lim_{h\to 0} \|\rho_h(f)\|=\|\widehat f\|_{\sup}$.\\
{\rm (b)} $\lim_{h\to 0} \|\rho_h(f)\rho_h(g)-\rho_h(f*g)\|=0$.
}

\Proof. We have 
$$\lim_{\hbar\to 0} \|\rho_{\hbar}(f)\| 
=\lim_{\hbar\to 0}\|\gvp_{\hbar}(\tilde f)\|=\|\gvp_{0}(\tilde f)\|
= \|\pi_0(f)\|=\|\widehat f\|_{\sup}
$$
and, for arbitrary $x$,
\begin{eqnarray*}
\|\rho_\hbar(f)\rho_\hbar(g)-\rho_\hbar(f*g)\|
&=& \|\pi_{(x,\hbar)}(\tilde f)\pi_{(x,\hbar)}(\tilde g)-\pi_{(x,\hbar)}(\widetilde{f*g})\|\\
&=&\|\pi_{(x,\hbar)}(\tilde f*\tilde g-\widetilde{f*g})\|\to \|\pi_0 (\tilde f*\tilde g-\widetilde{f*g})\|=0.
\end{eqnarray*}
\eproof

\section{Manifolds with Boundary}

\setcounter{equation}{0}

In the following, we shall denote by $X$ a smooth $n$-dimensional 
manifold with boundary, $\partial X$.
We assume that $X$ is embedded in a boundaryless manifold $\widetilde X$
and write $X^\circ$ for the interior of $X$. 
We also fix a Riemannian metric on $X$, so that we have $L^2$ spaces. 
We will show later on that the construction is independent of the choice of the metric. 
First of all, however, it is helpful to study the case where  $X=\R^n_+=\{ (x_1,\ldots, x_n)| x_n\geq 0 \}$ 
(including $x_n=0$!). We adopt the usual notation by writing an element $x \in \R_+^n$ as $x=(x',x_n)$.

\subsection*{Local Computation of the Asymptotic Green Term}
We change the formula for the $\hbar$-scaled boundary pseudodifferential operators with 
Fourier transformed symbol $f\in C^\infty_c(T\R_+^n)$ to
\begin{eqnarray}
\rho_\hbar(f)\xi(x)&=&\int_{x_n\geq \hbar v_n}f(x,v)\xi(x-\hbar v)dv\nonumber\\
&=&\hbar^{-n}\int_{w_n\ge0}f\big(x,\frac{x-w}\hbar\big)\xi(w)\,dw\label{hscaled}
,\quad \xi\in L^2(\R^n_+).
\end{eqnarray} A straightforward computation shows that 
\begin{eqnarray*}
(\rho_\hbar(f)\rho_\hbar (g)\xi )(x)&=& \int_{x_n\geq \hbar w_n}
\left( \int_{x_n\geq \hbar v_n}f(x,v)g(x-\hbar v,w-v)dv \right)\xi(x-\hbar w)dw\\
& =&\int_{x_n\geq \hbar w_n}\Big( \int f(x,v)g(x-\hbar v,w-v)dv \\ &&
  - \int_{x_n \leq \hbar v_n } f(x,v)g(x-\hbar v,w-v)dv \Big) \xi(x-\hbar w)dw,
\end{eqnarray*} 
where in the last line $f,g$ have to be understood as extended to functions in $C^\infty_c(T\R^n)$. The term 
\begin{eqnarray}\label{f*hg}
f*_hg= \int f(x,v)g(x-\hbar v,w-v)dv
\end{eqnarray}
is just the usual composition of Fourier transformed symbols of 
pseudodifferential operators on manifolds without boundary. 
We call the remainder, i.e. the operator which maps $\xi$ to 
\begin{eqnarray}
x&\mapsto&{-\int_{x_n\geq \hbar w_n}
\int_{x_n\leq \hbar v_n} f(x,v)g(x-\hbar v,w-v)dv\, \xi(x-\hbar w)dw}\nonumber\\
&=&-\int_{y_n\ge0}
\int_{x_n\leq \hbar v_n} f(x,v)g(x-\hbar v,y'-v',\frac{x_n}\hbar -y_n-v_n)\,dv\, \xi(x'-\hbar y',\hbar y_n)\, dy\mbox{\quad\quad }
\end{eqnarray} 
the  ``asymptotic Green'' term, because it corresponds to the leftover 
term in the composition of two truncated pseudodifferential operators 
in Boutet de Monvel's calculus, which is a singular Green operator, 
cf.\ \cite{MR53:11674}.   
In order to analyze it, we introduce the following notation: 

\begin{dfn}{kappa} 
{ For $0< \hbar \leq 1$ define 
$$\kappa_\hbar :C^\infty_c(T\R^{n-1}\times \R_+\times \R_+\times [0,1])\to\cL (L^2(\R^n_+))$$ 
 by
$$\kappa_\hbar (K)\xi (x)=\int_{y_n\geq 0} K(x',y',\frac{x_n}\hbar ,y_n,\hbar )\xi(x'-\hbar y',\hbar y_n)dv .$$
}
\end{dfn}

The asymptotic Green thus is of the form 
$\gk_\hbar(l_\hbar (f,g))$ with 
\begin{eqnarray}\label{lfg}
\lefteqn{l_\hbar (f,g)(x',y',x_n,y_n)}\nonumber\\
&=&
-\int_{x_n\leq v_n}f(x',\hbar x_n,v)g(x'-\hbar v',\hbar(x_n-v_n),y'-v',x_n-y_n-v_n)dv.
\end{eqnarray} 
As $\hbar\to0$ this tends to 
\begin{eqnarray}
 l(f,g)(x',y',x_n,y_n)=-\int_{x_n\leq v_n}f(x',0,y'-v',v_n)g(x',0,v',x_n-v_n-y_n)dv.
\end{eqnarray}
In fact, the difference $l_\hbar(f,g)-l(f,g)$ is an element of $C^\infty_c(T\R^{n-1}\times\R_+\times\R_+\times[0,1])$
which vanishes for $\hbar=0$. Similarly, the difference $f*_\hbar g-f*g\in C^\infty_c(T\R^n_+\times[0,1])$
vanishes for $\hbar=0$.   
\bigskip

In order to extend Theorem \ref{1.14} to manifolds with boundary, 
the asymptotic Green term has to be taken into account. 

\dfn{pi0}{For $\hbar =0$ we introduce
$$\pi_0^\partial :C^\infty_c(T\R^n_+\times [0,1])\oplus 
C^\infty_c(T\R^{n-1} \times \R_+\times \R_+\times [0,1])
\to \cL(L^2(T\R^{n-1}\times \R_+))$$ 
given by
\begin{eqnarray*}
\lefteqn{\pi_0^\partial (f\oplus K)\xi(x',v',v_n)=}\\
&=&\int_{w_n\geq 0} \left(f(x',0,v'-w',v_n-w_n,0) + K(x',v'-w',v_n,w_n,0)\right)\xi (x',w',w_n)dw
\end{eqnarray*} 
We complement $\pi_0^\partial$  by the map $\pi_0:C^\infty_c(T\R^n_+)\to \cL(L^2(T\R^n_+))$ 
in \eqref{tilderho}.
}

The crucial point is: 

\lemma{starprod}{The map
$$(\pi_0 ,\pi_0^\partial ):C^\infty_c(T\R^n_+)\oplus C^\infty_c(T\R^{n-1}\times \R_+\times \R_+)\to 
\cL (L^2(T\R_+^n)\oplus L^2(T\R_+^{n-1}\times \R_+))$$  
given by 
$$(\pi_0,\pi_0^\partial )(f\oplus K)=(\pi_0(f),\pi^\partial_0(f\oplus K))$$
turns $C^\infty_c(T\R_+^n)\oplus C^\infty_c(T\R^{n-1}_+\times \R_+\times \R_+)$ into an algebra. 
We denote this product with $*'$. Note that $f*'g=f*g+l(f,g)$.
}

It is clear that Theorem \ref{1.14}(b) will not remain true literally. Instead we obtain:

\thm{2.2}{For two symbols $f,g\in C^\infty_c(T\R^n_+)$ the following holds
$$\lim_{\hbar\to0} \| \rho_\hbar (f)\rho_\hbar (g)-\rho_\hbar(f*g)-\kappa_\hbar(l(f,g))
\|=0.$$
}

As in the case of manifolds without boundary, 
this will be related to 
the continuity of a  field of $C^*$-algebras which we will now introduce 

\begin{dfn}{field}{We denote by $A$ the  $C^*$-completion of 
$$A^\infty=C^\infty_c(T\R^n_+\times [0,1])\oplus C^\infty_c(T\R^{n-1}\times \R_+\times \R_+\times [0,1])$$
 in the representation $\rho_\hbar + \kappa_\hbar$, for $\hbar \not= 0$ and  $(\pi_0 , \pi_0^\partial )$ for $\hbar =0$,
i.e. under the mappings
$$f\oplus K\mapsto \left\{ \begin{array}{ll}
\rho_\hbar(f)+\gk_\hbar(K),&\hbar\not=0;\\\pi_0(f)\oplus \pi_0^\partial(f\oplus K)),&\hbar=0.
\end{array} \right. 
$$
There are obvious maps 
$$\varphi_\hbar :A  \to A(\hbar ),$$
where $A (\hbar )$ is the completion of $C^\infty_c(T\R^n_+)\oplus C^\infty_c(T\R^{n-1}\times \R_+\times \R_+)$  with respect to the specific representation in $\hbar$.}
\end{dfn}

We will show:

\thm{rncase}{The triple $(A, \{ A (\hbar ),\varphi_\hbar \}_{\hbar \in [0,1]})$ 
is a continuous field of $C^*$-algebras with $A(\hbar )$ 
isomorphic to the compact operators for $\hbar \not= 0$.
}

For fixed $\hbar \not=0$, the operators $\rho_\hbar(f)+\gk_\hbar(K)$
are compact, because they 
are integral operators with a square integrable kernel, 
 so  $A(\hbar)$
is isomorphic to the compact operators.  

We shall next analyze the field in more detail.
We abbreviate  
$$T=T\R^n_+\quad\text{ and }\quad \cT_\partial=T\R^{n-1}\times \R_+\times \R_+$$
and start with the following observation:

\begin{prop}{stjernealg}{As a subset of $A$, $A^\infty$ is a  
is a $*$-algebra.}
\end{prop}

\Proof. First we prove closure under multiplication. 
The product of $K_1,K_2 \in C^\infty_c(\cT_\partial\times [0,1])$ is  
just the convolution product of 
the two functions on the groupoid $\cT \R^{n-1}\times \R_+\times \R_+$, thus again a function in  
$C^\infty_c(\cT_\partial\times [0,1])$.

For  $f,g\in C^\infty_c({T})$ we  have already computed, cf.\ \eqref{f*hg}, \eqref{lfg}: 
$$\rho_\hbar (f)\rho_\hbar (g)=\rho_\hbar (\tilde{f}*_\hbar \tilde{g})+\kappa_\hbar (l_\hbar (\tilde{f},\tilde{g})).$$
where $\tilde{f},\tilde{g}$ are smooth extensions of $f,g$ to functions in
$C^\infty_c(T\R^n\times [0,1])$.
Since
$$( \pi_0 ,\pi_0^\partial)(f)(\pi_0,\pi_0^\partial)(g)=(\pi_0(f*g),\pi_0^\partial(f*g)+\pi_0^\partial(l(f,g) )$$
we see the closure under products of $f,g$. 

Checking the closure under products of $f$'s with $K$'s is straightforward.  
The same is true for the closure under involution.
\eproof

\subsection*{The Algebra in 0}

The algebra in zero, $A(0),$ is the completion of  
$$A(0)^\infty :=(C^\infty_c({T})\oplus
C^\infty_c(\cT_\partial),*')$$
in the representation $(\pi_0 , \pi_0^\partial )$. The summand 
$C^\infty_c(\cT_\partial)$ becomes an ideal in
$A(0)^\infty$.  We thus get the short exact sequence
\begin{equation} \label{glatkort}
0\to C^\infty_c(\cT_\partial) \to A(0)^\infty
\stackrel{q}\to C^\infty_c({T})\to 0. 
\end{equation}
As noted in the proof of Proposition \ref{stjernealg}, 
the algebra structure on  $C^\infty_c(\cT_\partial)$ 
comes from the groupoid structure on $\cT_\partial$, where  $\R_+\times \R_+$ 
carries the pair groupoid structure. Likewise, the algebra structure on
$C^\infty_c({T})$ stems from the
groupoid structure on ${T}$. Note that both groupoids are
amenable. 

\lemma{sesforA0}{We have a short exact sequence of $C^*$-algebras
\begin{equation}\label{stjernekort}
0\to C^*_r(\cT_\partial) \to A(0) \to
C^*_r({T}) \to 0.
\end{equation}}
\Proof. In the short exact sequence \eqref{glatkort}, the projection $q$, 
mapping  $f\oplus K$ to $f$, is a $*$-homomorphism. 
The trivial estimate 
$$\|\pi_0(f)\|_{\cL(L^2({T}))}\le 
\| \pi_0(f)\oplus \pi_0^\partial(f\oplus K)\|_{\cL(L^2({T})\oplus L^2(T\R^{n-1}\times\R_+))},$$
shows that $\pi$ extends
to a map $A(0)\to C^*_r({T})$ with  $C^*_r(\cT_\partial)$ in its kernel.
Since we may estimate the norm of $\pi_0^\partial(f)$ by the norm of $\pi_0(f)$,
we obtain \eqref{stjernekort}.\eproof\medskip

Alternatively, the lemma may be proven using only the amenability of the groupoids, similarly as in 
the proof of Theorem \ref{rn}, below.
Note that, via the Fourier transform, 
$$ C^*_r(\cT_\partial)\simeq
C_0(T^*\R^{n-1})\otimes \cK (L^2(\R_+))$$
and 
$$C^*_r({T})\simeq C_0(T^*\R_+^n).$$

\subsection*{Upper Semi-continuity}

\dfn{as}{On $A$ we define 
$$\| a\|_{as}=\max (\limsup_{\hbar \to 0}  \|\varphi_\hbar (a) \|,\|\varphi_0(a)\|).$$
This is a $C^*$-seminorm  which is continuous with respect to the norm of $A$. 
The quotient
$$A[0]=A/I, \quad \text{where}\quad
I=\{ a\in A |\ \| a\|_{as}=0 \},$$
therefore carries two norms: the quotient norm and $\|\cdot\|_{as}$. Both are equivalent by
\cite[Proposition 1.8.1]{Dixmier}, so that  $A[0]$ is a $C^*$-algebra with norm $\|\cdot \|_{as}$.
}
Since $\|a\|_{as}\ge \|\gvp_0(a)\|$ we have a natural map 
$$\Phi:A[0]\longrightarrow A(0).$$  

\begin{lemma}{vurdering}
 {Elements in  $A^\infty$
which are $0$ for $\hbar=0$ belong to $I$.}
\end{lemma}

\Proof. For 
$f\oplus K \in A^\infty $ 
it is easy to estimate
$$\|\rho_\hbar (f)\|\leq M_f\|f(\cdot , \hbar)\|_\infty \hbox{ and }
\| \kappa_\hbar (K) \| \leq M_K\|K(\cdot,\hbar )\|_\infty,$$
where $M_f$ and $M_K$ are constants depending on $f$ and $K$,  respectively, but not on $ \hbar$. 
\eproof

\thm{rn}{The field $(A,\{ A(\hbar ),\varphi_\hbar \}_{\hbar \in [0,1]})$ is upper semi-continuous in $0$.}
 
\Proof.
We denote by $R$ the closure of the range of the natural map $\gg:C_c^\infty(\cT_\partial)\to A[0]$.
This is an ideal in $A[0]$: Indeed, 
$C_c^\infty(\cT_\partial)$ is an ideal in $A(0)^\infty$, and 
the extension (e.g. constant in $\hbar$) of functions in $A(0)^\infty$ to functions in $A^\infty$ 
furnishes an embedding of $A(0)^\infty$ into $A[0]$ with dense range. 

Since $\cT_\partial$ is amenable, the quotient map $C^*(\cT_\partial)\to 
C^*_r(\cT_\partial)$ is an isomorphism.  It factorizes through 
$R$, since $R$ gives us a Hilbert space representation of $\cT_\partial$, 
while $\|a\|_{as}\ge\|\gvp_0(a)\|$.

This leads to a commutative diagram of natural maps
\begin{eqnarray*}
&&C^*(\cT_\partial)\\
&\nearrow&~~~\downarrow~~~\\
C_c^\infty(\cT_\partial) & \stackrel{}{\hookrightarrow} &~~~R\subseteq A[0], \\
& \searrow &~~~ \downarrow~~\\
&&C_r^*(\cT_\partial)
\end{eqnarray*}
where the upper vertical arrow is surjective, since the inclusion has dense range.
The invertibility of the quotient map implies that the lower  vertical arrow is an
isomorphism.

Next we define a map $\tilde{q}:A[0]\to C^*_r(T)$: By definition, $A[0]$ is the set of 
equivalence classes of Cauchy sequences in $A^\infty$ with respect to $\|\cdot\|_{as}$. 
Given such a Cauchy sequence $a_k=(f_k\oplus K_k)$, we may evaluate at $\hbar=0$ and obtain 
a sequence $(f_k^0\oplus K_k^0)$ in $A(0)^\infty$. As $\|a_k\|_{as}\ge \|\gvp_0(a_k)\|$,
the sequence $(f^0_k)$ is a Cauchy sequence in $C_r^*(T)$; moreover, the mapping 
$(a_k)\mapsto (f^0_k)$ is well-defined and continuous.
In view of Lemma \ref{vurdering} its kernel is $R$.   
 
Combining this with the short exact sequence \eqref{stjernekort}
we obtain the following commutative diagram of short exact sequences
\begin{equation}
\begin{array}{ccccccccc}\label{cd}
0 & \to & C_r^*(\cT_\partial) & \to & A[0] & \stackrel{\tilde{q}}\to &C_r^*({T})&\to &0 \\
 \| &&\| &&\downarrow \Phi&&\|&&\| \\
0 & \to & C_r^*(\cT_\partial) & \to & A(0) & \to &C_r^*({T})&\to &0
\end{array}.
\end{equation}
We conclude from the five lemma that $\Phi$ is an isomorphism, and therefore 
$$\limsup_{\hbar \rightarrow 0} \|\varphi_\hbar (a)\|\leq \|\varphi_0(a)\|,$$
i.e. the field is upper semi-continuous in $0$.\eproof\medskip

What is still missing is the proof of the lower semi-continuity of the field $A$. 
It will be given at the end of Section 2, since there is no simplification 
for the half-space case.

\subsection*{The Tangent Groupoid for a Manifold with Boundary}

\dfn{2.1}{We denote by $T^-X$ the subset of  
$T\widetilde X$ formed by all vectors  $(m,v)\in T\widetilde X|_{X}$ 
for which $\exp_m(-\gve v)\in X$ for sufficiently small $\gve>0$.  
This is a semi-groupoid with addition of vectors.
Note that $T^-X=TX^\circ\cup T^-X|_{\partial X}$ 

We define $\cT^-X$ as the disjoint union $T^-X \cup (X\times X\times ]0,1])$, 
endowed with the fiberwise semi-groupoid structure 
induced by the semi-groupoid structure on $T^-X$ and the groupoid structure on
$X\times X$.
As in the boundaryless case, 
we glue $T^-X$ to $X\times X \times ]0,1]$ via the charts
$$T^-X\times [0,1]\supseteq U\ni (m,v,\hbar)\mapsto \left\{ 
\begin{array}{cc}
(m,v) & \hbox{for }\hbar =0 \\
(m, \exp_m(-\hbar v),\hbar)& \hbox{for }\hbar \not= 0
\end{array}
\right.$$
and
let $\cT^-X(0)=T^-X$ and $\cT^-X(\hbar)=X\times X\times \{\hbar\}$.

In order to avoid problems with the topology of $\cT^-X$ (which is in general not a manifold with corners)
we let $C^\infty_c(\cT^-X)=C^\infty_c(\cT\widetilde{X})|_{\cT^-X}$.

}

\subsection*{C*-algebras Associated to the Semi-groupoids $T^-X$ and $\cT^- X$}
We start with $T^-X$. 
Let  $ \Ctc ( T^-X)$ denote the smooth functions on $T^-X$ which have compact support in $T^-X$.
In analogy with Definition \ref{kappa} we introduce  
\begin{eqnarray*}
\pi_0&:&\Ctc(T^-X)\to \cL(L^2(T X^\circ ))\quad\text{and}\\
\pi_0^\partial&:&\Ctc(T^-X)\to \cL(L^2(T^-X|_{\partial X}))
\end{eqnarray*} 
acting by 
\begin{eqnarray}
\pi_0(f)\xi(m,v)=\int_{T_mX}f(m,v-w)\xi(m,w)\, dw,\\
\pi_0^\partial(f)\xi (m,v)=\int_{T^+_mX}f(m,v-w)\xi(m,w)\, dw.\label{kappaf}
\end{eqnarray} 
Note that due to its compact support in $T^-X$, the function $f$ 
naturally extends (by zero) to $TX$.
\dfn{Cr}{We denote by $C^*_r(T^-X)$ the $C^*$-algebra generated by 
$\pi_0$ and $\pi_0^\partial$, i.e.~by the map
$\Ctc(T^-X)\ni f\mapsto 
(\pi_0(f),\pi_0^\partial(f))\in \cL(L^2(TX^\circ)\oplus L^2(T^-X|_{\partial X}))$.}

At first glance, this definition seems to overlook the operators of the form
$\pi_0^\partial(K)$ in \ref{kappa} and operators of the form $\pi_0(f)$ and $\pi_0^\partial (f)$, where $f\in C^\infty_c(T\tilde{X})|_{TX}$. In fact, this is not the case. 
The second type of operators belongs to $C^*_r(T^-X)$, because we take the closure under the adjoint operation and addition. 
The reason that the first type of operators is in $C^*_r(T^-X)$, is the 
well-known relation between operators of half-convolution and 
Toeplitz operators, which we recall, below. We denote by $\fT$ the algebra of all 
Toeplitz operators on $L^2(S^1)$ and by $\fT_0$ the ideal of all operators whose 
symbol vanishes in $-1$.

\lemma{Toeplitz2}{Let $f\in C^\infty_c(\R).$ Then the operator 
$$L^2(\R_+)\ni\xi\mapsto \left(s\mapsto \int_0^\infty f(s-w)\xi(w)\,dw\right)\in L^2(\R_+)$$
is unitarily equivalent to the Toeplitz operator $T_\gvp$ with symbol $\gvp(z)=\hat f(i(z-1)/(z+1)).$ 
Note that $\gvp(-1)=\hat f(\infty)=0$. 

The $C^*$-algebra generated by the
operators in the image of $C^\infty_c(\R)$ under this map 
is precisely the ideal $\fT_0$, while the compact operators in $\fT$  are generated 
by their commutators.}

\Proof. Plancherel's theorem shows that the above operator of half convolution is the truncated pseudodifferential
operator with symbol $\hat f$, mapping $\xi\in L^2(\R_+)$ to $\op(\hat f)_+\xi(s)=
\int e^{ist}\hat f(t)\widehat{(e^+\xi)}(t)\,dt$, where $e^+\xi$ is the extension (by zero) of $\xi$ 
to $\R$. 

Now one observes that the unitary $U:L^2(S^1)\to L^2(\R_+)$ given by 
$Ug(t)=\frac{\sqrt 2}{1+it}~g\left(\frac{1-it}{1+it}\right)$ maps the Hardy space $H^2$ to 
$F(L^2(\R_+))$ with the Fourier transform $F$, and that $\op(\hat f)_+$
is $F^{-1}UT_\gvp U^{-1}F$. See \cite[Section 2]{RS} for details. 

For the second statement, one first notes that the 
$C^*$-algebra generated by these operators is a subalgebra of 
$\fT_0$. On the other hand, $\fT_0$ consists of the operators of the form 
$T_\gvp+C$, where $\gvp\in C(S^1)$ vanishes in $-1$, and $C$ is compact.
According to \cite[Proposition 7.12]{Douglas}, the commutators of all $T_\gvp$, $ \gvp\in C(S^1)$, 
generate the compacts, hence so do the commutators of those $T_\gvp$, where 
$\gvp$ vanishes in $-1$. As these $T_\gvp$ can be approximated by elements in the image of 
$C^\infty_c(\R)$, the proof is complete. \eproof

\lemma{2.12}{We have a representation $\pi_0^\partial$ of 
$C^\infty_c(T\partial X\times \R_+\times \R_+)$ on $L^2(T^-X|_{\partial X})$ via 
\begin{eqnarray}\label{kappaK}
\pi_0^\partial(K)\xi(m,v',v_n)=\int K(m,v'-w',v_n,w_n) \xi(m,w',w_n)\,dw'dw_n.
\end{eqnarray} 

The closure of its range is isomorphic to 
$$J=C_0 (T^*\partial X)\otimes \cK (L^2(\R_+)).$$

$J$ is an ideal in $C_r^*(T^-X)$ generated by commutators of elements of the form $\pi_0^\partial(f)$. }

\Proof. The algebraic tensor product 
$C^\infty_c(T\partial X)\otimes C^\infty_c(\R_+\times \R_+)$ 
is dense in $C^\infty_c(T\partial X\times \R_+\times \R_+)$. 
Due to the continuity of 
$$\pi_0^\partial :C^\infty_c(T\partial X\times \R_+\times \R_+) \to \cL(L^2(T^-X|_{\partial X})$$
it is sufficient to determine the closure of $\pi_0^\partial ( C^\infty_c(T\partial X)\otimes C^\infty_c(\R_+\times \R_+))$.

It is clear that  $\pi_0^\partial ( C^\infty_c(T\partial X)\otimes C^\infty_c(\R_+\times \R_+))\subseteq J$. 
In fact, we have equality, 
since the Fourier transform gives an isomorphism $C_r(T\partial X)\to C_0(T^*\partial X)$ 
and since 
a compact operator on $L^2(\R_+)$ can be approximated by a 
Hilbert-Schmidt operator, thus by an integral operator with kernel in $C^\infty_c(\R_+\times \R_+)$.

In order to see that $J$ is contained in $C^*_r(T^-X)$, 
it is  sufficient to approximate both factors of a pure tensor $h\otimes c$, where
$h\in C_0(T^*\partial X)$ and $c\in \cK(L^2(\R_+))$. 
For the first task we choose a function in $C^\infty_c(T\partial X)$ whose 
fiberwise Fourier transform is close to $h$ in sup-norm.  
For the second, we  refer to Lemma \ref{Toeplitz2}. 
In particular, we see that $J$ also is generated by commutators.

A direct computation shows that $J$  is an ideal in $C^*_r(T^-X)$. 
\eproof\medskip

\dfn{kappa_h}{We let
$$\cut =
C^\infty_c(TX)\oplus C^\infty_c(T \partial X\times \R_+\times\R_+).$$ 
This is a dense $*$-subalgebra of $C^*_r( T^-X)$. We will denote the product in this subalgebra by $*'$.
\\

For $\hbar \not= 0$ we obtain representations of 
$C^\infty_c(\cT^-X)=C^\infty_c(\cT \widetilde{X})|_{\cT^-X}$ in $\cL(L^2(X))$ by:
\begin{equation} \label{repsg}
 \pi_\hbar (f)\xi (m)=\frac{1}{\hbar^n}\int
f(m,\tilde m,\hbar)\xi (\tilde m)d\tilde m.
\end{equation}
Note that these are the natural groupoid representations for $X\times X\times ]0,1]$.

We denote by $C_r^*(\cT^-X)$  the reduced $C^*$-algebra generated by  
$\pi_\hbar$, $0\le \hbar\le 1$, and $\pi_0^\partial$. }

For $X=\R^n_+$ we have $\cut =A(0)^\infty$,   $C^*_r(T^-X)=A(0)$ and $C^*_r(\cT^-X)=A$. 
Also there are evaluation maps 
$$\varphi_\hbar :C^*_r(\cT^- X)\to C_r^*(\cT^- X)(\hbar ). $$ 

\thm{hoved}{We have \begin{eqnarray*}
C^*_r(\cT^-X(\hbar ))&=&\cK (L^2(X)),\quad \hbar \not= 0;\\
C^*_r(\cT^- X(0))&=&C^*_r(T^-X).\end{eqnarray*} 
 Moreover:
$( C^*_r(\cT^- X),\{ C^*_r(\cT^- X(\hbar )),\varphi_\hbar\}_{\hbar \in [0,1]})$ is a continuous field of $C^*$-algebras.
}

The first two statements are obvious. 
For the proof of upper semi-continuity, we will essentially follow the ideas for the half-space case. 
Our first task is the construction of a representation of $\cut$.
To this end, we will simply extend $f\in C^\infty_c(TX)$ and $K\in C^\infty_c(T\partial X\times\R_+\times \R_+)$ to functions 
$\tilde f$ and $ \tilde K$ on $\cT^-X$ as described below, then apply \eqref{repsg}.

Choose  a function $\psi \in C^\infty (X\times X)$
which is one on a neighborhood of the diagonal, $0\leq \psi \leq 1$, such that 
$$\exp :T^-X \rightarrow X\times X$$ 
maps a neighborhood of the zero section diffeomorphically to the support of $\psi$.  

For $f\in C^\infty_c(TX)$ we define  $\tilde{f} \in C^\infty_c(\cT^-X)$ by 
\begin{equation} \label{udvidf} \tilde{f}(m,\tilde m,\hbar)=\psi (m,\tilde m) 
f \left( m,  -\frac{\exp^{-1}(m,\tilde m)}{\hbar } \right).
\end{equation}

We next identify a neighborhood $U$ of $\partial X$ in $X$ with $\partial X\times [0,1[$
and write $U\ni m =(m',m_n)$ with $m'\in\partial X $ and $m_n\ge 0$. 
We also  choose a function $\chi \in C^\infty_c(X)$ supported in $U$ with  $0 \leq \chi
\leq 1$ and $\chi\equiv 1$ near $\partial X$. 
For $K \in C^\infty_c(T\partial X \times \R_+\times\R_+ )$ we then
define $\tilde{K} \in C^\infty_c(X\times X\times]0,1])$ by
\begin{eqnarray}\label{udvidk} 
\tilde K(m,\tilde m,\hbar)=
\chi(m)\chi(\tilde m)\psi(m,\tilde m)~K\left(m',-\frac{ \exp^{-1}(m',\tilde m')}\hbar
,\frac{m_n}{\hbar},\frac{\tilde m_n}{\hbar }\right). 
\end{eqnarray}

\rem{halfspace}{In the half-space case with the flat metric 
we have, for fixed $f$ and $K$,
$$\pi_\hbar(\tilde f)=\rho_\hbar(f)\quad\text{and}\quad \pi_\hbar (\widetilde{K})  
=\kappa_\hbar(K)$$  
provided $\hbar$ is sufficiently small.} 

\cor{Relation_to_(0.1)}{We then obtain the analog of Property \eqref{0.1}: 
\begin{eqnarray}
\lim_{\hbar\to 0} \|\pi_\hbar(\widetilde f)+\pi_\hbar(\widetilde K) \| =\max\{\|\pi_0(f)\|,\|\pi_0^\partial(f\oplus K)\|\}.
\end{eqnarray}
}

\extra{metric}{Metrics}{The construction of $C^*_r(\cT^-X)$ and the extensions  (\ref{udvidf}), (\ref{udvidk}) used a metric, 
but  $C^*_r(\cT^-X)$ is independent of the choice: 
Let $\nu_1, \nu_2$ be two different metrics on $X$, 
and denote by $\mu_1,\mu_2$ the associated measures on $X$ 
as well as the fiberwise measures in $TX$. 
Let $k \in C^\infty  (X)$ be given by $$\mu_1=k\mu_2.$$
Multiplication by $\sqrt{k}$ yields a unitary 
$$U: (L^2(X),\mu_1)\rightarrow (L^2(X),\mu_2),$$
and multiplication by $\sqrt{k(m)}$ a family of unitaries 
$$U_m :(L^2(T_m^-X),\mu_1) \rightarrow (L^2(T_m^-X),\mu_2).$$
We define 
$$\phi:C^\infty_c(\cT^-X)\to C^\infty_c(\cT^-X)$$
taking $f(m,v,0)$  to $f(m,v,0)k(m)$ for $\hbar=0$ and $f(m,\tilde m,\hbar)$ to 
$f(m,\tilde m,\hbar)\sqrt{k(m)k(\tilde m)}$, $\hbar\not=0$.
Then $\pi^1_\hbar(f)=U^{-1}\pi_\hbar^2(\phi(f))U$, 
where $\pi^{1}_\hbar$ and $\pi_\hbar^2$ are the representations
induced by $\mu_1$ and $\mu_2$. A corresponding relation holds for  $\pi_0^\partial$. 
Hence  $C^*_r(\cT^-X)$ is independent of the metric. 
}
The following lemma clarifies the influence of the extension by different metrics.

\begin{lemma}{metrik}{
Let $f\in C^\infty_c(TX)$. 
Denote by $\tilde{f}^i$ the extension of $f$ with respect to the metric $\nu_i$, $i=1,2$. 
Then
$$\| \pi_\hbar (\phi (\tilde{f}^1)) -\pi_\hbar (\widetilde{\phi (f)}^2)\| \rightarrow 0 \hbox{ for } \hbar \rightarrow 0.$$ 
Here $\pi_\hbar$ is understood with respect to $\mu_2$.}
\end{lemma}

\Proof. This follows from Lemma \ref{vurdering}, 
since $\phi (\tilde{f}^1)-\widetilde{\phi (f)}^2$ is a function in $C^\infty_c(\cT^-X)$ 
which is zero at $\hbar =0$.
\eproof
\medskip

A similar statement holds if we start with $K \in C^\infty_c(T\partial X \times \R_+\times\R_+ )$. 

\subsection*{Upper Semi-continuity}
We again use the seminorm 
$$\| a\|_{as}=\max \{ \| \varphi_0(a)\| ,\limsup_{\hbar \to 0}\| \varphi_\hbar (a)\|\}$$ 
for elements in $C^*_r(\cT^-X)$ and introduce the analog of $A[0]$: 
$$C^*_?(T^-X)=C^*_r (\cT^- X)/I,$$
where

$$I=\{ a \in C^*_r( \cT^- X)|~ \| a\|_{as} =0 \}.$$ 

The notation $C^*_? (TX^- )$ is justified by the following:

\begin{prop}{fundamental}{The mappings $f\mapsto\tilde f$ and $K\mapsto\tilde K$ induce a 
$*$-homomorphism $\Psi$ from $(C^\infty_c(T^-X),*')$ to $C^*_?(T^-X)$ with dense range,
and we have 
\begin{eqnarray}\label{AsMult}
 \lim_{\hbar \to 0}\|\pi_\hbar(\tilde f)\pi_\hbar(\tilde g)-\pi_\hbar(\widetilde{f*'g)}\|=0,\quad f,g\in C^\infty_c(TX).
\end{eqnarray}
}
\end{prop}

\Proof.  Choose an open covering $\{ U_i\}$ of $X$, 
where each $U_i$ can be identified with an open subset of $\R^n$ or $\R^n_+$.
By possibly shrinking the $U_i$, we may assume that 
the function $\psi$
used in \eqref{udvidf} and \eqref{udvidk} equals $1$ on $U_i\times U_i$ and that
the function $\chi$ is $\equiv 1$ on $U_i$  whenever
$U_i$ intersects the boundary.
We also fix a subordinate partition of unity $\{\psi_i\}\subset C^\infty_c(U_i)$.

For $f,g\in C^\infty_c(TX)$ we have 
$(\psi_if)*'g=(\psi_if)*'(\eta_ig)$ for
each $\eta_i\in C^\infty_c(U_i)$ with $\psi_i\eta_i=\psi_i$. Moreover,  
$\pi_ \hbar(\widetilde {\psi_if})\pi_\hbar (\tilde g)=
\pi_ \hbar(\widetilde {\psi_if})\pi_\hbar (\widetilde{\gt_i g})$ 
for suitable $\gt_i\in C^\infty_c(U_i)$, provided $\hbar $ is small.
Hence 
\begin{eqnarray}\label{normi}
\lefteqn{ \| \pi_\hbar ( \widetilde{f*'g} )-
\pi_\hbar (\tilde{f}) \pi_\hbar (\tilde{g})\|
\leq \sum\| \left(\pi_\hbar  ( \widetilde{(\psi_if)*'g} )
-\pi_\hbar (\widetilde{\psi_if}) \pi_\hbar (\tilde{g})\right)\|}\nonumber\\
&=&\sum \|  \pi_\hbar  (\widetilde{\psi_if*'\eta_ig}))
-\pi_\hbar (\widetilde{\psi_if}) \pi_\hbar (\widetilde{\gt_ig})\|.
\end{eqnarray}   
For sufficiently small $\hbar$, all operators will have support in $U_i\times U_i\times [0,1]$
so that we are working on Euclidean space.
According to Lemma \ref{metrik} we can also,  modulo terms converging to zero as $\hbar \rightarrow 0$, 
use the Euclidean metric.
So we are precisely in the situation considered at the beginning of the section. 
The explicit computation shows that 
\begin{eqnarray}\label{diff}
\pi_\hbar(\widetilde{f*'g}))-\pi_\hbar(\tilde f)\pi_\hbar (\tilde g)=\rho_\hbar(f*g-f*_\hbar g)+\gk_\hbar(l(f,g)-l_\hbar(f,g)).
\end{eqnarray} 
As $f*g-f*_\hbar g\in C^\infty_c(T\R^n_+\times[0,1])$ and 
$l(f,g)-l_\hbar(f,g)\in C^\infty(T\R^{n-1}\times\R_+\times\R_+\times[0,1])$
vanish for $\hbar=0$,  
the difference \eqref{diff} is in $I$ by Lemma \ref{vurdering}.
Hence \eqref{normi} tends to zero, and $\Psi (f*' g)=\Psi (f)\Psi (g)$.  
The remaining $*$-algebra properties are checked similarly.

In order  to see that the image of $\Psi$ is dense in $C^*_? (T^-X)$,
we simply note that the evaluation at $\hbar=0$ associates to an element $F$ in $C^\infty_c(\cT^-X)$
an element in $\cut$ whose extension via \eqref{udvidf}, \eqref{udvidk} 
induces the same element in $C_?^*(\cT^-X)$ 
by  Lemma \ref{vurdering}.\eproof

\rem{AM}{Property \eqref{AsMult} is the analog of the asymptotic 
multiplicativity \eqref{0.2} in the case of manifolds with boundary. 
In particular, we have established Theorem \ref{2.2}.  }

With Proposition \ref{fundamental}, 
the proof of the following theorem is analogous to 
that of Theorem \ref{rn}.
\thm{uppergeo}{
$( C^*_r(\cT^- X),\{ C^*_r(\cT^- X)(\hbar ),\varphi_\hbar\}_{\hbar \in [0,1]})$ 
is  upper semi-continuous in $0$.}

\forget{\extra{AM}{Asymptotic multiplicativity}{Proposition \ref{fundamental} in connection with 
Lemma \ref{vurdering} shows that
\begin{eqnarray*}
 \lim_{\hbar \to 0}\|\rho_\hbar(f)\rho_\hbar(g)-\rho_\hbar({f*'g)}\|=0,\quad f,g\in C^\infty_c(\cT^-X),
\end{eqnarray*}
where $f*'g$ denotes the product of $f$ and $g$ in $C_r^*(\cT^-X)$.

For $f,g\in C^*_r(\cT^-X)$ we choose sequences $(f_k), (g_k)$ in $C^\infty_c(T^-X).$
Upper semi-continuity then establishes asymptotic multiplicativity for $f,g\in C^*_r(\cT^-X)$.}
}

\subsection*{Lower Semi-continuity}

As in the classical case \cite{LandsmanRamazan} 
lower semi-continuity is proven by introducing strongly continuous representations
using the groupoid structure. 
We split the representations into two: 
One taking care of the contribution from the interior of the manifold, 
i.e. the convolution part, and one taking care of the boundary part, 
i.e. half convolution and kernels on the boundary.

For the lemmata, below, we note that -- by construction -- 
$\pi_0$ and $\pi_0^\partial$ extend to $C^*_r(\cT^-X)$.

\begin{lemma}{lavet}{$\liminf_{\hbar \to 0}\| \varphi_\hbar ( a)\|\geq \| \pi_0 (a)\|$  for all $a\in  C^*_r(\cT^-X)$. }
\end{lemma}

\Proof. 
According to Proposition \ref{fundamental} it is sufficient to show that 
\forget{, the families of the form $\rho_\hbar (\tilde f+\tilde{K})$
form a dense subset of 

It follows from Proposition \ref{fundamental} that the families of the form $\rho_\hbar (f+\tilde{K})$, 
$\hbar \not= 0$ and $(\pi_0(f),\pi_0^\partial(f)+\pi_0^\partial(K))$, where $f\in C_c^\infty (\cT^-X)$, 
and $\tilde{K}$ is constructed from $K\in C^\infty_c(\cT\partial X \times \R_+\times\R_+)$ as in 
(\ref{udvidk}), form a dense subset of $C^*_r(\cT^-X)$. 
Indeed, the class of $a\in C^*_r(\cT^-X)$ in $C^*_? (T^-X)$ 
can be approximated by elements of the form $\rho_\hbar (f)+\rho_\hbar(\tilde{K})$. 
By an analog of Lemma \ref{vurdering}, elements in $I$ can be approximated by elements in 
$C^\infty_c(TX\times [0,1])$ which vanish at $\hbar=0$. Hence it sufffices to show}
\begin{equation} \label{in}
\lim_{\hbar \to 0}\| \rho_\hbar (\tilde f+\tilde{K} )\| \geq \| \pi_0 (f)\| \quad\text{for}~~f\oplus K\in \cut.
\end{equation}

For  $g \in C^\infty_c (\cT^- X)$   define  
$$\| g \|_{\infty ,\hbar}^2 =\sup_{m\in X} \left\{ \frac{1}{\hbar^n}\int_X | g(x,m,\hbar )|^2 dx \right\} \hbox{ for } \hbar \not= 0,$$
and  
$$\|g \|^2_{\infty , 0} = \sup_{m \in X} \left\{ \int_{T_mX}|g(m,v,0 )|^2dv\right\} ,  \hbox{ for } \hbar=0.$$ 
Set
$$\|g\|_\infty =\sup_{\hbar \in [0,1]}\|g \|_{\infty , \hbar}.$$
It is easily checked that
\begin{eqnarray}\label{norm1h}
\| \pi_\hbar (\tilde f+\tilde{K})\|
=\sup \Big\{ \Big\|\frac 1 {\hbar^n} 
\int  (\tilde f( \cdot ,m,\hbar)+\tilde{K}(\cdot,m,\hbar ))g( m,\cdot,\hbar )\,dm\Big\|_{\infty ,\hbar } 
~\Big|~ \|g\|_\infty \leq 1\Big\}
\end{eqnarray} 
for $\hbar \not= 0$, and 
\begin{eqnarray}\label{norm0}
 \|\pi_0 (f)\| 
= \sup \Big\{ \Big\|  \int f(\cdot,v,0)g(\cdot ,\cdot-v,0)dv\Big\|_{\infty ,0} 
~\Big|~ \|g\|_\infty \leq 1 \Big\}:
\end{eqnarray}
In fact, for \eqref{norm1h} we note that ``$\ge$'' follows from the estimate
\begin{eqnarray*}
\lefteqn{\left\|\frac1{\hbar^n}\int\tilde f(m_1,m,\hbar)g(m,m_2,\hbar)\, dm\right\|_{\infty,\hbar}^2 =
\left\|\frac1{\hbar^n}\pi_\hbar(\tilde f)g(\cdot,m_2,\hbar)\right\|_{\infty,\hbar}^2}\\
&=&\sup_{m_2\in X}\frac1{\hbar^n} 
\left\|\pi_\hbar(\tilde f)g(\cdot,m_2,\hbar)\right\|_{L^2(X)}^2
\le \left\|\pi_\hbar(\tilde f)\right\|^2 \sup_{m_2\in X}\frac1{\hbar^n} 
\left\|g(\cdot,m_2,\hbar)\right\|_{L^2(X)}^2\\
&=&\left\|\pi_\hbar(\tilde f)\right\|^2\|g\|_{\infty,\hbar}^2\le
\left\|\pi_\hbar(\tilde f)\right\|^2\|g\|_{\infty}^2 .
\end{eqnarray*} 
For the reverse inequality we choose $g(x,m,\hbar)=s(m)\xi(x)\hbar^n\gvp(\hbar)$, 
where $s\in C^\infty_c (X)$, $s\le 1$,
$\|\xi\|_{L^2(X)}=1$ with  $\|\pi_\hbar(\tilde f)\xi\|\ge \|\pi_\hbar(\tilde f)\|-\gve$, and  
$\gvp\in C^\infty_c(]0,1])$ is equal to one outside a neighborhood of zero. 
Equation \eqref{norm0} follows by a similar argument.
 
Now suppose that  $g\in C^\infty_c(\cT^-X)$ and $g(x,m,h)=0$ for $x\in\partial X$. 
Then the weak convergence of $\tilde{K}$ towards zero implies that
\begin{eqnarray*}
\lim_{\hbar \rightarrow 0} \Big\| \frac 1 {\hbar^n} \int  (\tilde f( \cdot ,m,\hbar )+
\tilde{K}(\cdot, m,\hbar ))g( m , \cdot,\hbar )\,dm\Big\|_{\infty , \hbar}
&=& \|  \int f(\cdot , v,0)g(\cdot ,\cdot-v,0)dv\|_{\infty ,0}.
\end{eqnarray*}
As the set of these $g$ is dense in $\{g\in C^\infty_c(\cT^-X)\ |\ \|g\|_\infty\le 1\}$,  \eqref{in} follows.
\eproof

\begin{lemma}{lavto}{$\liminf_{\hbar \to 0} \| \varphi_\hbar (a)\|
\geq \| \pi_0^\partial(a)\|$ for all $a\in C^*_r(\cT^-X)$.}
\end{lemma}

\Proof. As in the proof of Lemma \ref{lavet} we only have to show that
\begin{equation} \label{into}
\liminf_{\hbar \to 0} \| \rho_\hbar (\tilde f+\tilde{K})\| \geq \| \pi_0^\partial (f\oplus K)\|,
\end{equation}
for $f \in C^\infty_c(TX)$ and $K\in C^\infty_c(T\partial X\times \R_+\times\R_+ )$. 

We let $P_\hbar$ be the projection in $L^2(X)$ given by multiplication by the characteristic 
function of $\partial X\times[0,a_\hbar[$, where 
$$a_\hbar \rightarrow 0 \hbox{ for } \hbar \rightarrow 0
\text{~ and ~} \frac{a_\hbar} \hbar \rightarrow \infty \hbox{ for } \hbar \rightarrow 0.
$$ 
As $\|P_\hbar \pi_\hbar (\tilde{f}+\tilde{K}) P_\hbar \| \leq \| \pi_\hbar(\tilde{f}+\tilde{K}) \| $,
it is enough to show that  
$$\liminf_{\hbar \rightarrow 0}
\| P_\hbar \pi_\hbar(\tilde f+\tilde{K}) P_\hbar \|\geq  \| \pi_0^\partial(f\oplus K)\|.$$
Since we are free too choose a metric, we fix a metric on $\partial X$ and the standard metric on $[ 0,a_\hbar [$. 

As in the proof of Lemma \ref{lavet}, 
we equip the space $C^\infty_c(\cT \partial X \times [0,\infty [)$ with norms 
\mbox{$\| \cdot \|_{\infty ,\hbar}$}, \mbox{$\|\cdot\|_\infty$}, 
which are just like the norms before, on $\cT \partial X$ instead of $\cT^- X$, combined with the $L^2$-norm on $[0,\infty[$. 
For $f\in C^\infty_c(\cT^-X)$ and $K\in C^\infty_c (T \partial X \times \R_+\times\R_+)$
we define representations on   $C^\infty_c(\cT \partial X \times [0,\infty [)$ by
\begin{eqnarray*}
 \eta_\hbar (f)g(m_1,m_2,\hbar, b)&=&\frac 1 {\hbar^{n-1}}\int_{a\in [0,\frac{a_\hbar}\hbar]} 
f(m_1,\hbar b,m,\hbar a,\hbar)g( m, m_2 , \hbar, a )dmda,\ \ b\in[0,\frac{a_\hbar}{\hbar}[,\hbar\not=0; \\
\eta_0 (f)g(m_1,v,0,b)&=&\int_{T_{m_1}\partial X\times \R_+}f(m_1,0,v-w,b-a,0)g(m_1,w,0,a)\,dwda;\\ 
\eta_0 (K)(m_1,v,0,b)&=&\int_{T_{m_1}\partial X\times \R_+} K(m_1,v-w,b,a)g(m_1,w,0,a)\, dwda.
\end{eqnarray*} 
Note that 
$\|P_\hbar \pi_\hbar(f) P_\hbar\|=\|D_{\hbar} P_\hbar \pi_\hbar (f)P_\hbar D_{\hbar^{-1}}\|=\sup \{ \|\eta_\hbar (f)g\|_{\infty ,\hbar}~|~\|g\|_\infty \leq 1 \}$, 
where $D_\hbar$ is the dilation operator in the normal direction, given by $D_\hbar f(x',x_n)=f(x',\hbar x_n)$.

As before
$$ \| \pi_0^\partial (f\oplus K)\| =\sup \{ \|(\eta_0 (f\oplus K))g\|_{\infty,0} |\|g\|_\infty \leq 1\}.$$

Plugging in the definitions of $\tilde{f}$ and $\tilde{K}$ (omitting the cut off functions) we get
$$\eta_\hbar (\tilde{f})g(m_1,m_2,\hbar, b)=\frac{1}{\hbar^{n-1}}\int_{[0,\frac{a_\hbar }{\hbar}]}f\left( m_1,\hbar b,-\frac{\exp^{-1} (m_1,m)}{\hbar} ,b-a\right) g(m,m_2,\hbar,a) dmda$$$$
\eta_\hbar (\tilde{K})g(m_1,m_2,\hbar,b)=\frac{1}{\hbar^{n-1}}\int_{[0,\frac{a_\hbar }{\hbar}]}K\left( m_1, b,-\frac{\exp^{-1} (m_1,m)}{\hbar} ,a\right) g(m,m_2,\hbar,a)dmda\\
$$
Using dominated convergence and the fact that $g$ for small $\hbar$ looks like $g_0(m,-\frac{\exp^{-1} (m,m_2)}{\hbar} ,\hbar ,a)$, $g_0\in C^\infty_c(T\partial X \times [0,1]\times [0,\infty [)$, 
we get 
$$\lim_{\hbar \to 0} 
\| (\eta_\hbar (\tilde f+\tilde{K})g\|_{\infty,\hbar}=\|(\eta_0 (f\oplus K)g\|_{\infty,0},$$
and  (\ref{into}) follows.
\eproof
\\ 

Lemma \ref{lavet} and \ref{lavto} imply that 
$\liminf_{\hbar \to 0} \| \varphi_\hbar (a)\|\geq \|\varphi_0 (a)\|,$
i.e.

\thm{nedad}{
$( C^*_r(\cT X^-),\{ C^*_r(\cT X^-)(\hbar ),\varphi_\hbar\}_{\hbar \in [0,1]})$ 
is  lower semi-continuous in $0$.}

This finishes the proof of Theorem \ref{hoved}.

\section{$K$-theory of the Symbol Algebra $C^*_r(T^-X)$}
\setcounter{equation}{0}
 $C_c^\infty (TX^\circ)$ with the fiberwise convolution product is a $*$-ideal of $\cut$. 
After completion,  $C^*_r(TX^\circ )$ becomes a $C^*$-ideal of $C^*_r(T^-X)$, and we have a short exact 
sequence  
\begin{equation} \label{foelge}
0 \to  C^*_r(TX^\circ) \to C^*_r(T^-X) \to C^*_r(T^-X)/ C^*_r(TX^\circ)\to 0.
\end{equation}

\prop{Toeplitz}{%
The quotient $Q= C^*_r(T^-X)/ C_r^*(TX^\circ)$ is naturally isomorphic to 
$C_0(T^*\partial X)\otimes \fT_0$ for the ideal $\fT_0$ of the 
Toeplitz algebra introduced before Lemma \ref{Toeplitz2}.
}
 
\Proof. Define $$\Psi:C_c^\infty(T^-X)\to \cL(L^2(T^-X|_{\partial X}))\text{~ by ~}\Psi(f\oplus K)=\pi_0^\partial(f)+\pi_0^\partial(K)$$
with the maps in \eqref{kappaf} and \eqref{kappaK}. This is a $*$-homomorphism with respect to $*'$,
and  $C^\infty_c(TX^\circ)$ is in its kernel.
We first show that $\ker\Psi=C^*_r(TX^\circ)$:
Since $C_r^*(T^-X)$ is the closure of $C^\infty_c(T^-X)$ with respect to the norm 
$$\|f\oplus K\|=\max\{|\pi_0(f)\|,\|\pi_0^\partial(f)+\pi_0^\partial(K)\|\},$$
and $C_r^*(TX^\circ)$ is the closure of $C^\infty_c(TX^\circ)$ with respect to $\|\pi_0(f)\|$, 
we have $C^*_r(TX^\circ)\subseteq \ker\Psi$. 

On the other hand, suppose that $a\in\ker\Psi$; i.e., $a$ is the equivalence class of a Cauchy sequence
$(f_k\oplus K_k)\in C^\infty_c(T^-X)$ with $\pi_0^\partial(f_k)+\pi_0^\partial(K_k)\to 0$.
We next note that 
\begin{eqnarray*}
\|\pi_0(f_k)\|&=&\sup\{|\hat f_k(m,\gs)|~|~(m,\gs)\in T^*X\}:
\text{~and~}\\
\|\pi_0^\partial(f_k)+\pi_0^\partial(K_k)\|&\ge& \sup\{|\hat f_k(m,\gs)|~|~(m,\gs)\in T^*X|_{\partial X}\}
\end{eqnarray*}
Indeed the first inequality follows from the fact that, via fiberwise Fourier transform, 
$\pi_0(f_k)$  is equivalent to multiplication by $\hat f_k(m,\sigma)$.
For the second, we observe first that $\|\pi_0^\partial(f_k)\|=\sup\{|\hat f_k|\}$ as a consequence 
of the fact that translation of $\xi=\xi(m,w)$ in the direction of $w_n$ preserves 
$\|\pi_0^\partial(f_k)\xi\|$ in $L^2(T^-X|_{\partial X})$.
On the other hand, $\pi_0^\partial(K_k)\xi=0$ provided we translate sufficiently far. Hence
$\|\pi_0^\partial(f_k)+\pi_0^\partial(K_k)\|\ge \|\pi_0^\partial(f_k)\|$. 

We conclude that the fiberwise Fourier transforms $\hat f_k$ tend to zero uniformly on $T^*X|_{\partial X}$.
Hence the Cauchy sequence $(f_k)$ may be replaced by an 
equivalent Cauchy sequence $(g_k)$ with $g_k\in C^\infty_c(TX^\circ)$.
We conclude that $\pi_0^\partial(K_k)\to 0$ so that $(K_k)\sim 0$, and therefore
$\ker\Psi\subseteq C^*_r(TX^\circ)$.
 
Hence $\Psi$ descends to an injective $C^*$-morphism on $Q$; in particular, it has closed range.

Now we observe that we have a natural identification of 
$TX|_{\partial X}$ with $T\partial X\times \R$ and consequently 
of $T^-X|_{\partial X}$ with $T\partial X\times \R_-$. 
Hence $\cL(L^2(T^-X|_{\partial X}))\cong \cL(L^2(T\partial X)\otimes L^2(\R_-)).$ 

Suppose that, at the boundary,  $f\in C^\infty_c(TX)$ is of the form $f(x',0,v',v_n)= g(x',v')h(v_n)$
with $g\in C^\infty_c(T\partial X)$ and $h\in C_c^\infty(\R)$. Then $\pi_0^\partial(f)=\pi_0^{\partial,0} (g)\otimes\pi_0^{\partial,n} (h)$,
where $\pi_0^{\partial,0} $ is the convolution operator by $g$, acting on $L^2(T\partial X)$, while  
$\pi_0^{\partial,n} (h)$ is the operator of half convolution acting on $L^2(\R_-)$ 
(note that $\R_-\cong T^-\R_+|_{\{0\}}$).
Via Fourier transform, the operator  $\pi_0^{\partial,0} (g)$ is unitarily equivalent to 
multiplication by $\hat g\in C_0(T^*X)$, while, 
according to Lemma \ref{Toeplitz2}, $\pi_0^{\partial,n} (h)$ is unitarily equivalent to a 
Toeplitz operator in $\fT_0$. The closure of the image of the span of the pure 
tensors thus gives us $C_0(T^*\partial X)\otimes \fT_0$. 

We know already from Lemma \ref{2.12} that -- via the Fourier transform --  
the image of $C^\infty_c(T\partial X\times\R_+\times\R_+)$ 
can also  be identified with a subset of  
$C_0(T^*\partial X)\otimes \cK\subseteq C_0(T^*\partial X)\otimes \fT_0$. 
This completes the argument.
\eproof

\thm{K}{Via fiberwise Fourier transform $C^*_r(TX^\circ )$ can be identified with 
$C_0(T^*X^\circ)$ and 
the inclusion $ C_0(T^*X^\circ) \cong C^*_r(TX^\circ )\hookrightarrow C^*_r(T^-X)$ induces an isomorphism of K-groups 
$$K_i ( C^*_r(T^-X))\cong K_i ( C_0(T^*X^\circ)),\quad i=0,1.$$}

\Proof. It is well-known (or easily checked) that $K_i(\fT_0 )=0$, $i=0,1$.  
Thus it follows from the Künneth formula that $K_i(C_0(T^*\partial X )\otimes \fT_0)=0$, $i=0,1$. 
The result now is a consequence of \eqref{foelge} and the associated six term exact sequence.
\eproof

\noindent Johannes Aastrup, Institut für Mathematik, Universität Hannover, Welfengarten 1, \mbox{30167 Hannover}, Germany,
email: {\tt aastrup@math.uni-hannover.de}\bigskip

\noindent Ryszard Nest, Department of Mathematics, Copenhagen University, Universitetsparken 5, \mbox{2100 Copenhagen}, Denmark,
email: {\tt rnest@math.ku.dk}\bigskip

\noindent Elmar Schrohe, Institut für Mathematik, Universität Hannover, Welfengarten 1, \mbox{30167 Hannover}, Germany,
email: {\tt schrohe@math.uni-hannover.de}
\end{document}